\documentclass[final,3p,times]{elsarticle} 

\usepackage{amssymb}

\usepackage{multirow}
\usepackage{algorithmic}
\usepackage{amsmath}
\usepackage{amsfonts}
\usepackage{amsthm}
\usepackage[english]{babel} 
\usepackage{dsfont}
\usepackage{amssymb}
\usepackage{graphicx}        
\usepackage{subcaption}
\usepackage{epstopdf}
\usepackage{epsfig}
\usepackage[table]{xcolor}
\usepackage{subfloat}
\usepackage{float}
\usepackage[thinlines]{easytable}
\usepackage{array}
\usepackage{mathtools}
\usepackage{hyperref}
\usepackage{comment}
\usepackage{bm}
\usepackage[shortlabels]{enumitem}
\usepackage{epstopdf}
\usepackage{booktabs}
\usepackage{url}
\usepackage[font=small,labelfont=bf,textfont=it,indention=.1cm,width=1.1\textwidth]{caption}
\usepackage[linesnumbered,algoruled,boxed,lined]{algorithm2e}

\newenvironment{myfont}{\fontfamily{lmtt}\selectfont}{\par}
\DeclareTextFontCommand{\textmyfont}{\myfont}

\SetKwComment{Comment}{/* }{ */}
\numberwithin{equation}{section}
\newtheorem{remark}{Remark}[section]

\newtheorem{lemma}{Lemma}[section]
\newtheorem{theorem}{Theorem}[section]

\newtheorem{proposition}{Proposition}[section]
\newtheorem{assumption}{Assumption}[section]

\newcommand{\Rd}{ {\mathbb{R}^\textup{d}}}
\newcommand{\Rm}{ {\mathbb{R}^\textup{m}}}
\newcommand{\R}{ \mathbb{R}}

\newcommand{\ve}{\varepsilon}
\renewcommand{\d}{\textup{d}}
\newcommand{\m}{\textup{m}}
\newcommand{\supp}{\textup{supp}}
\newcommand{\law}{\textup{Law}}

\newcommand{\refup}{\textup{ref}}

\newcommand{\argmin}{\textup{argmin}}
\newcommand{\diam}{\textup{diam}}
\newcommand{\V}{\textup{V}}
\newcommand{\proj}{\textup{proj}}

\begin{document}

\begin{frontmatter}



\title{Model predictive control strategies using
 consensus-based optimization\tnoteref{t1}}
\tnotetext[t1]{The authors thank the Deutsche Forschungsgemeinschaft (DFG, German Research Foundation) for the financial support under Germany’s Excellence Strategy EXC-2023 Internet of Production 390621612 and under the Excellence Strategy of the Federal Government and the Länder,  333849990/GRK2379
(IRTG Hierarchical and Hybrid Approaches in Modern Inverse Problems), 320021702/GRK2326,
 442047500/SFB1481 within the projects B04, B05  and B06 and through SPP 2410 Hyperbolic Balance Laws in Fluid Mechanics: Complexity, Scales, Randomness (CoScaRa) within the Project(s) HE5386/26-1 and HE5386/27-1.
Support through DFG HE5386/30-1 and  EU DATAHYKING No. 101072546 is also acknowledged.  The work of G.B. is funded by the Deutsche Forschungsgemeinschaft (DFG, German Research Foundation) through 320021702/GRK2326 ``Energy, Entropy, and Dissipative Dynamics (EDDy)''.}

 \author[aachen]{Giacomo Borghi\corref{cor1}}
 \ead{borghi@eddy.rwth-aachen.de}
 \author[aachen]{Michael Herty}\ead{herty@igpm.rwth-aachen.de}

 \cortext[cor1]{Corresponding author}

  \affiliation[aachen]{organization={RWTH Aachen University, Institut für Geometrie und Praktische Mathematik},
               addressline={Templergraben 55},
                           city={Aachen},
                                 postcode={52062},
             country={Germany}}

\begin{abstract}
Model predictive control strategies require to solve in an sequential manner, many, possibly non-convex,  optimization problems. In this work, we propose an interacting  stochastic agent system to solve  those problems.  The agents  evolve  in pseudo-time and in parallel  to the time-discrete state evolution. The  method is  suitable for non-convex, non-differentiable objective functions. The convergence properties are investigated through mean-field approximation of the time-discrete system, showing  convergence in the case of additive linear control. We validate the proposed strategy by applying it to the control of a  stirred-tank reactor non-linear system.
\end{abstract}




\end{frontmatter}


\section{Introduction}

Consensus-Based Optimization (CBO) is a class of stochastic multi-agent methods for the optimization of non-convex and possibly non-differentiable objective functions introduced in \cite{pinnau2017consensus}. Similar to  meta-heuristics algorithms such as Particle-Swarm Optimization (PSO) \cite{kennedy1995particle} or Genetic Algorithms (GA) \cite{tang1996genetic}, to name a few, CBO methods are based on a set of interacting agents (sometimes also called \textit{particles}) that explore the search space to determine  a global solution. The exploration behavior is achieved by inserting stochasticity in the dynamics, while exploitation of the known objective function landscape is achieved by a simple, derivative-free interaction between the agents. In particular,  the agents evolve according to a consensus-type interaction, leading to concentration at the agents attaining low values of the objective function. CBO algorithms benefit from rigorous mathematical frameworks that allow for an analytical understanding of the optimization dynamics and its convergence properties \cite{carrillo2018analytical, fornasier2021consensusbased}. This is typically achieved by performing an approximation of mean-field of the multi-agents system, see e.g. \cite{huang2021meanfield}.

So far, most of the CBO and related algorithms have been limited to finite-dimensional controls. In order to extend the applicability towards a setting of possibly time-continuous controls, we employ a Model Predictive Control (MPC) strategy based on piecewise constant controls. The MPC technique has been widely analysed in engineering and mathematical papers and we refer to  \cite{mayne2000constrained,camacho2013model,grune2017nonlinear,azmi} for a non-exhaustive list of references. 
In those works, typically, the time windows is subdivided. In each sub-interval, an open-loop control  is applied, after which the optimization procedure is repeated using the updated state.  Here, we couple the MPC strategy with a CBO optimization routine to solve a possibly non-linear control system in each time-interval of a  time-discrete settings. The agent dynamics evolves in parallel with the state evolution, but possibly on  a different (faster) time scale. The control applied to the state could then given by, for instance, a weighted  mean of the state of the CBO agents. This allows to use a single CBO multi-agent system to dynamically solve the possibly high-dimensional, sequential MPC sub-problems even having a  reduction of the computational cost.


After performing a mean-field approximation of the CBO dynamics, we theoretically analyze the convergence properties of the algorithm under the assumption of linear additive control. In these settings, we are able to prove the convergence of the agents' weighted mean towards the optimal control for all MPC subproblems. We are also able to estimate a maximum number of CBO iterations needed, for each sub-problem. 

In MPC, particle filters are typically used for the estimation of the state evolution in presence of uncertainties \cite{andrieu2004particle,botchu2007nonlinear}. When used for the control selection, a common strategy consists of treating the optimization problem as an inference problem solved with Monte Carlo methods \cite{doucet2002MCMC}, see, for instance, \cite{stahl2011} where a sequential strategy is also proposed. This approach is sometimes combined with the Simulated Annealing (SA) algorithm \cite{kirkpatrick1983annealing}, see \cite{kantas2009SMC,devilliers2011particle}. Differently to such particle-based methods, in the CBO dynamics the agents interact which each other  to solve the optimization problem. The CBO dynamics resembles the PSO  (applied for MPC problems in \cite{xu2016mpcpso}), but  it is  amendable for analytical considerations and shows computationally better performance \cite{borghi2023chaos}.


The paper is organized as follows. In Section \ref{sec:problem} we present the problem settings and recall the MPC strategy. Section \ref{sec:CBO} is devoted to the presentation of the novel CBO algorithm for the MPC problem, that is analyzed in Section \ref{sec:analysis} by mean-field approximation. To numerically validate the proposed strategy we illustrate numerical examples in Section \ref{sec:num}. Conclusive remarks and possible outlook are discussed in the last section.

\section{Problem settings}
\label{sec:problem}

We consider model predictive control for the discrete time, possibly non-linear, control system
\begin{equation}
x_{(n+1)} = \Phi(x_{(n)}, u_{(n)})
\label{eq:state}
\end{equation}
where $x_{(n)} \in \mathcal{X} \subseteq\Rm$ is the state, and $u_{(n)} \in \Rd$ the control. We call $\Phi: \mathcal{X} \times \Rd \to \mathcal{X}$ the \textit{transition} map. 
The (constant) control $u_n$ is assumed to be taken, component-wise, among a certain range of values $u_{\min}, u_{\max} \in \Rd, u_{\min} <u_{\max}$, and we denote the admissible control set by
\begin{equation}
\mathcal{U} = \{ u \in \Rd\,:\, u_{\min} \leq u \leq u_{\max}\}\,.
\end{equation}
Given that most mathematical models of real-life processes are continuous-in-time models, $\eqref{eq:state}$ can be intended as a sampled data system, that is, a model for the data collected from an underlying process at times $t_n = n \Delta t$ for some sampling period $\Delta t$. 

To determine the control strategy, we assume the data $\{x_{(n)}\}_{n\geq0}$ available by equation \eqref{eq:state} and assume the transition map $\Phi$ to be available for computations.  Let $\overline{n}$ be a given time horizon, and $x_{(0:\overline n)}  = \{x_{(n)} \}_{n=0}^{\overline n}$ $u_{(0:\overline n)}  = \{u_{(n)} \}_{n=0}^{\overline n}$ be the set of states and controls up to (time) $\overline n$. 
Given a  desired reference state $x^{\refup}:[0,\overline {n}\Delta t]\to \mathcal{X}$, the aim of the control is to minimize a loss function $L$  given by 
\begin{equation}
\begin{split}
\min_{u_{(0:\overline n)} \in \mathcal{U}^{\overline n}}\, L(u_{(0:\overline n)}):= \sum_{n = 1}^{\overline n+1 } |x_{(n)} -  x^{\refup}(t_{n})|^2 + \nu \sum_{n = 0}^{\overline n } |u_{(n)}|^2 \quad \textup{subject to} \quad \eqref{eq:state}\,.
\end{split}
\label{pb:total}
\end{equation}
The loss $L$ measures  the discrepancy between $x_{(n)}$ and $x^{\refup}(t_{n+1})$ and  we  consider a regularization term for the control, dependent on a parameter $\nu>0$.   

Minimizing \eqref{pb:total} requires to solve the optimization problem on the set of controls $u_{(0:\overline{n})}$. This is typically computationally expensive, as it requires to consider the  time interval $[0, \overline{n} \Delta t ]$ and to apply the forward operator $\overline{n}$ times for each evaluation of $L$.  Instead of solving \eqref{pb:total}, we use a Model Predictive Control (MPC) strategy. 

Let $n_p \ll \overline{n}$ be the prediction horizon, in MPC, the original optimization problem is divided in simpler sub-problems for the controls given by
\begin{equation}
\min_{u_{(n: n + n_p)} \in \mathcal{U}^{\overline n}}\, L_n (u_{(n:n + n_p)}):= \sum_{m = n+1}^{n + n_p +1 } |x_{(m)} -  x^{\refup}(t_{m})|^2 + \nu \sum_{m = n}^{n+n_p} |u_{(m)}|^2 \quad \textup{subject to} \quad \eqref{eq:state}\,,
\label{pb:mpc_0}
\end{equation}
which are solved sequentially for $n = 0, \dots, \overline{n}$. This is possible as $x_{(n)}$ only depends on $x_{(n-1)}$ and $u_{(n-1)}$. Once a solution $\tilde u_{n:n+n_p}$ is obtained,  we apply the first component $u_{(n)} = \tilde u_{(n)}$ to obtain state $x_{(n+1)}$ and re-iterate the procedure. Other choices are also possible, but we focus here on this setting.

%
%


\section{Consensus-Based Optimization  for MPC}
\label{sec:CBO}

As MPC method, we couple a Consensus-Based Optimization (CBO)  dynamics \eqref{eq:CBO} to  the state dynamics \eqref{eq:state}.   CBO algorithms are stochastic multi-agent systems that have been proved to be suitable for non-convex global optimization problems. During the algorithm iteration, a set of $N \in \mathbb{N}$ agents, explore the search space according to a nonlinear interaction rule inspired by consensus dynamics in social interactions.
 
In the following, $n$ will be the index of time associated to the state dynamics, while $k$ indicates the CBO  index of the pseudo-time step. We will show that  $n$ and $k$ can be  coupled. Consider  a fixed time  $n$ and state variable $x_{(n)}$. The associated objective function $L_n$ of the  MPC  problem is given by \eqref{pb:mpc_0}.  
The CBO agents of problem $n$ are represented by a set of random variables $\{U_{(n, k)}^i\}_{i=1}^N$ taking values in $\R^{\d  n_p}$. 
At every time $(n,k)$ the CBO agents move towards a weighted average, the consensus point, given by
\begin{equation}\label{eq:malphaN}
m_{L_n}^\alpha[f^N_{(n,k)}] = \frac{\sum_{i=1}^N U^i_{(n,k)} \exp \left(-\alpha L_n(U_{(n,k)}^i)\right)}{\sum_{i=1}^N \exp \left(-\alpha L_n(U_{(n,k)}^i)\right)},
\end{equation}
where $f^N_{(n,k)} = (1/N)\sum_{i=1}^N \delta_{U_{(n,k)}^i}$ is the empirical distribution associated with the agent system at time $(n,k)$. 
Note that for $\alpha>0$, the exponential weights are larger for agents with low objective values, and correspond to the Boltzmann-Gibbs distribution of  the objective function $L_n$.

To ensure the agents are constrained to the feasible set $\mathcal U$, we  apply the projection operator component--wise 
\begin{equation*}
\proj_{\mathcal{U}}(v):= \underset{u \in \mathcal{U}}{\argmin}\;|u - v|.
\end{equation*}
Finally, the  update rule for each agent  is given by the  scheme 
\begin{equation}
\begin{cases}
U_{(n,k+1/2)}^i = U_{(n,k)}^i + \lambda \tau \left (m^\alpha_{L_n}[f^{N}_{(n,k)}] - U_{(n,k)}^i \right)  + \sigma \sqrt{\tau} D \odot \theta_{(n,k)}^i \\
U_{(n,k+1)}^i = \proj_{\mathcal{U}}\left(U_{(n,k+1/2)}^i \right)
\label{eq:CBO}
\end{cases}
\end{equation}
starting by (arbitrary) initial data $\{U_{(n,0)}^i\}_{i=1}^N$. Constants $\lambda, \sigma, \tau>0$ are parameters of the algorithm to be defined later, and $\theta_{(n,k)}^i \sim \mathcal{N}(0,I_{\d n_p})$ is sampled form the standard normal distribution.  

Dynamics \eqref{eq:CBO} is characterized by a deterministic step towards the consensus point $m_{L_n}^\alpha[f_{(n,k)}^N]$, and a random component whose variance depends on $D\in \R^{\d n_p}_{\geq 0}$, where $D$ may depend, in turn, on the difference $m_{L_n}^\alpha[f^{N}_{(n,k)}] - U_{(n,k)}^i$. We  also consider non-degenerate diffusion operators $D$ for which $D>0$ component-wise , see e.g. \cite{carrillo2022sampling}. 


Let $\overline{k}$ be the time horizon of the CBO computation. We consider the computed optimal control to be given by the weighted mean $m_{L_n}^\alpha[f_{(n,\overline k)}]$. In particular,
let $\tilde u_{(n:n+n_p)} = m_{L_n}^\alpha[f_{(n,\overline k)}]$ be the computed solution to \eqref{pb:mpc_0}. Then, we  apply the control $u_{(n)} = \tilde u_{(n)}$ in the dynamics. Afterwards, we continue with  the next control problem.

Assuming stability of the state dynamics as well as of the reference state $x^\refup$, we initialize the CBO  system for the optimal control problem $n$ with the final state configuration of problem $n-1$. This corresponds to consider equation \eqref{eq:CBO} and initial conditions
\begin{equation}
U_{(n, 0)}^i = U^i_{(n-1, \overline{k})} \qquad \textup{for all}\;\; i = 1, \dots, N\,.
\label{eq:CBOinitial}
\end{equation}
The resulting MPC-CBO method is described in Algorithm \ref{alg:mpc-cbo}.

 \begin{algorithm}
\caption{Model Predict Control strategy coupled with Consensus-Based Optimization routine.}
\label{alg:mpc-cbo}

\textbf{Input:} $x_{(0)}, \Phi, x^\refup, \Delta t,\nu, n_p\;$ \;
Set CBO parameters $\lambda, \sigma>0, \lambda\tau\in (0,1), N, \overline{k} \in \mathbb{N}\;$\;
Initialize CBO  system $\{U_{(0,0)}^i \}_{i=1}^N \subset \R^{\d n_p}\;$\;
\For{$n = 0, 1, \dots, \overline{n}-1$}{
\For{$k = 0, 1, \dots, \overline{k}$}{
  Update  $U_{(n, \overline{k})}^i$ according to \eqref{eq:CBO} for all $i = 1, \dots, N$\;
  }
  Compute  $m_{L_n}^\alpha[f_{(n, \overline{k})}]$ according to \eqref{eq:malphaN}\;
  $\tilde u_{(n:n+n_p)} = m^\alpha_{L_n}[f_{(n, \overline{k})}]\;$\;
  $u_{(n)} = \tilde u_{(n)}\;$\;
  $x_{(n+1)} =  \Phi(x_{(n)},u_{(n)})$ \;
   $U_{(n+1, 0)}^i = U^i_{(n, \overline{k})}\;$ for all $i = 1, \dots, N$\;
  }
\end{algorithm}

Some remarks are in order.


\begin{remark} \label{r:timecontinuous}
Note that for each single time step $n\to n+1$ of the evolution of the state $x_{(n)}$ , the CBO agent system evolves for $\overline{k}$ steps. In the continuous setting this corresponds to two different time scales. 

\par
 
The state update \eqref{eq:state} is determined by  sampling or simulation procedure of an underlying time-continuous system $x_{(n)} = x_{t_n}, t_n = \Delta n$, with $dx_t/dt = f(x_t, u_t)$. In this context, it interesting to verify if Algorithm \ref{alg:mpc-cbo} admits a time-continuous description.
If we keep the control horizon $T_p  = n_p\Delta t$ fixed,  the MPC loss at time $t$ formally reads 
\[ 
L_{t}(u_{[t,t+T_p]}) = \int_{t}^{t+T_p} \left (|x_s - x^\refup|^2 + \nu |u_s|^2\right) ds
\]
for $\Delta t \to 0$,
leading formally to an optimization problem over functions $u_{[0,T_p]}$.

\par

The CBO update has been considered in the literature as an Euler-Maruyama discretization \cite{platen1999intro} of a system of $N$ It{\^o}'s Stochastic Differential Equations (SDEs) for a time-continuous $\{U_s^i\}_{i=1}^N$ process with time step $\tau>0$ \cite{pinnau2017consensus}. The proposed method can also be formulated in such a framework (without projection step, for simplicity):
let $\Delta t/ \tau = \bar k \in \mathbb{N}$ be fixed, consistently with the choice  \eqref{eq:CBOinitial} we introduce the times
\[ t_{n,k} = n\Delta t +   k \tau\,, \]
and set  $U_{t_{n,k}}^i = U_{(n,k)}^i$.
By formally taking the limit $\Delta t, \tau \to 0$ with $\Delta t/\tau = \overline{k},$ we obtain the system
 \begin{equation}
 \begin{cases}
dx_{t}/dt = f\left(x_t,u_t = m_{L_t}^\alpha[f_t^N](0)\right) \\
dU_t^i (s)= \frac 1\epsilon \left(\lambda (m_{L_t}^\alpha[f_t^N](s) - U^i_t(s))dt + \sigma D\odot dB_t^i (s)\right) \quad \textup{for all}\; s\in[0,T_p]
\end{cases}
 \end{equation}
where, following the standard notation of two-scale dynamics, we introduced $\epsilon =1/\overline{k}$, and where $B_t^i$ is a family of parametrized Brownian processes.

\end{remark}

\section{Convergence analysis}
\label{sec:analysis}


In this section, we investigate the convergence properties of Algorithm \ref{alg:mpc-cbo} by relying on a mean-field approximation of the CBO agent dynamics. For the theoretical analysis, we assume  the following simplifying assumptions.
\begin{assumption} 
\label{asm:linearu} We consider a prediction horizon $n_p = 1$ and assume that the transition map $\Phi$ is linear with respect to the control, i.e., 
\begin{equation}
\Phi(x,u) = \Phi_s(x) + F_c u  
\end{equation}
for  $\Phi_s: \mathcal{X} \to \mathcal{X}$ and a matrix $F_c \in \R^{\d\times\m}$. We  assume that $F_c$ has 
$F_c^\top F_c + \nu I_\d$
 full-rank. As for the noise in the CBO agent dynamics, we assume that  
$D = \mathbf{1}_\d$. 
\end{assumption}

The CBO dynamics read  then
\begin{equation}
U_{(n,k+1)}^i = \proj_{\mathcal{U}} \left( U_{(n,k)}^i + \lambda \tau \left (m^\alpha_{L_n}[f^{N}_{(n,k)}] - U_{(n,k)}^i \right) + \sigma \sqrt{\tau} \theta_{(n,k)}^i\,\right).
\label{eq:CBO_nondeg}
\end{equation}

We derive the mean-field approximation for $N \gg 1$. Then, we establish convergence guarantees for single step $n$ of the MPC strategy. Finally, we show independence with respect to $n$ and extend the convergence for all problems $n = 0, \dots, \overline{n}-1$ considered.

\subsection{Mean-field approximation}

To study the large time evolution of the CBO agent system \eqref{eq:CBO_nondeg}, we assume propagation of chaos \cite{sznitman1991chaos} of the marginals $F^N_{(n,k)} \approx (f_{(n,k)})^{\otimes N}$ for $N \gg1$, where $f_{(n,k)} \in \mathcal{P}(\Rd)$. Under such assumption, the consensus point satisfies \cite{pinnau2017consensus}
\begin{align*}
m^\alpha_{"_n}[f^N_{(n,k)}]   
 \approx 
 m_{L_n}^\alpha[f_{(n,k)}]
\end{align*}
where we extended the definition of $m_{L_n}^\alpha[\cdot]$ for a probability measure $f \in \mathcal{P}(\Rd)$, as
 \begin{equation}
m_{L_n}^\alpha[f]:= \frac{\int u \exp(-\alpha L_n(u))f(du)}{\int \exp(-\alpha L_n(u))f(du)}\,.
\label{eq:malpha}
\end{equation}
We note that update  \eqref{eq:CBO_nondeg} is symmetric with respect to  the particular index $i$
\begin{equation}
\overline{U}_{(n,k+1)} = \proj_{\mathcal{U}} \left(\overline{U}_{(n,k)} + \lambda \tau\left( m^\alpha_{L_n}[f_{(n,k)}] - \overline{U}_{(n,k)}\right) + \sigma \sqrt{\tau}\,\overline{ \theta}_{(n,k)} \right),
\label{eq:discr_mono}
\end{equation}
where $\overline{\theta}_{(n,k)} \sim \mathcal{N}(0, I_\d)$ and $f_{(n,k)} = \textup{Law}(\overline{U}_{(n,k)})$.

For notational simplicity, we neglect  the iterative step $k$ and introduce the following mean-field dynamics
\begin{equation}
\begin{split}
u' &= \proj_{\mathcal{U}}\left( u + \lambda \tau \left(m^\alpha_{\tilde L}[f] - u\right) + \sigma \sqrt{\tau} \theta \right) \\
& =: \mathcal{C}^{CBO}(u,f, \theta, \tilde L)
\end{split}
\label{eq:collCBO}
\end{equation}
for a given $f\in \mathcal{P}(\Rd)$ and $u,\theta\in \Rd$. For any measurable test function $\phi: \Rd \to \R$,  the mono-agent process satisfies 
\[
\mathbb{E} \left [\phi (\overline{U}_{(n,k+1)}) \right] = \mathbb{E} \left [\phi (\mathcal{C}^{CBO}(\overline{U}_{(n,k)}, f_{(n,k)}, \theta, L_{(n)}) \right] 
= \mathbb{E} \left [\phi (\overline{U}_{(n,k)}') \right] 
\]
or, equivalently, the time-discrete mean-field evolution
\begin{equation}
\int \phi(u)\, f_{(n,k+1)}(du) = \iint \phi(u')\, f_{(n,k)}(du)\, d\mu(\theta)
\label{eq:MFdiscr}
\end{equation}
where, here, $u' = \mathcal{C}^{CBO}(u,f_{(n,k)},\theta, L_{n})$ is given by \eqref{eq:collCBO}. For the sake of readability, we will favor the notation $u'$ over the more precise $\mathcal{C}^{CBO}(u,f_{(n,k)},\theta,L_n)$ , whenever the intention  is clear.

\begin{remark} For quantitative estimates of the mean-field approximation error introduced we refer to \cite{gerber2023mean}. If we assume the initial data are sampled from a common distribution, the expected error has been shown in similar constrained settings to be of order $\mathbb{E}|m^{\alpha}_{L_n}[f_{(n,k)}^N] - m^{\alpha}_{L_n}[f_{(n,k)}]| = \mathcal{O}(N^{-1})$ \cite{fornasier2020hypersurfaces}.
\end{remark}

\subsection{Convergence for each single MPC sub-problem}

In the following, we consider a fixed $n\in \mathbb{N}_0 := \mathbb{N} \cup \{ 0\}$ and determine conditions such that the  CBO dynamics converges towards the global solution of the MPC sub-problem $n$.

Under Assumption \ref{asm:linearu}, sub-problem \eqref{pb:mpc_0} reduces to the constrained convex quadratic programming problem
\begin{equation} 
\min_{u \in \mathcal{U}}\; L_n (u):= \big| \Phi_s (x_{(n)}) + F_c u -  x^{\refup}(t_{n+1}) \big |^2 + \nu | u|^2\,.
\label{pb:mpc}
\end{equation}

\begin{lemma} Let Assumption \ref{asm:linearu} hold, then there exists a unique solution $u^*_{(n)}$. Let $(\eta_{(n)}^{*,1},\eta_{(n)}^{*,2}) \in \Rd\times \Rd$ be the correspondent Lagrange multipliers, and let $\lambda_{A}^{\textup{min}}$, $\lambda_{A}^{\textup{max}}>0$ be the smallest and largest eigenvalue of $A: = F_c^\top F_c + \nu I_\d$. 

For all $u \in \mathcal{U}$, it holds
\begin{equation}
\label{eq:grow2}
\lambda_{A}^{\min} | u - u_{(n)}^*|^2 \leq L_n(u)  - L_n(u_{(n)}^*)  \leq \left(|\eta_{(n)}^{*,1}|+|\eta_{(n)}^{*,2}| \right)|u - u_{(n)}| +  \lambda_{A}^{\max} | u - u_{(n)}^*|^2\,.
\end{equation}
\end{lemma}

\begin{proof} Solution uniqueness follows from the convexity of the problem. The first-order necessary and sufficient conditions are 
\begin{equation} \label{eq:kkt}
2\left( F_c^\top F_c + \nu I_\d \right) u_{(n)}^* - 2F_c^\top \left( \Phi_s (x_{(n)})  -  x^{\refup}(t_{n+1})  \right) + \eta_{(n)}^{*,1} - \eta_{(n)}^{*,2} = 0\,
\end{equation} 
with $\eta_{(n)}^{*,1}, \eta_{(n)}^{*,2}\in \Rd_{\geq0}$ being the Lagrange multipliers for the box constraints $u -  u_{\max} \leq 0, u_{\min} - u \leq 0$, respectively.

Let $\tilde x := \Phi_s(x_{(n)}) - x^\refup(t_{n+1})$. We note that second order Taylor expansion at $u_{(n)}^*$  leads to
\begin{align*}
L_n(u) - L_n(u_{(n)}^*) &= 2\left\langle (F^\top_cF_c + \nu I_\d)u^*_{(n)} - F_c^\top \tilde x, u - u^*_{(n)}   \right  \rangle + \left \langle u - u_{(n)}^*, \left(F^\top_cF_c + \nu I_\d \right)(u - u_{(n)}^*) \right \rangle \\
& = 2\left\langle - \eta_{(n)}^{*,1} + \eta_{(n)}^{*,2}, u - u^*_{(n)}   \right  \rangle + \left \langle u - u_{(n)}^*, \left(F^\top_cF_c + \nu I_\d \right)(u - u_{(n)}^*) \right \rangle
\,.
\end{align*}
The upper bound follows from Cauchy-Schwarz inequality. For the lower bound, we note that \begin{align*}
- \langle\eta^{1,*}, u - u_{(n)}^* \rangle &= - \langle\eta^{1,*}, u - u_{\max} \rangle  - \langle\eta^{1,*}, u_{\max} - u_{(n)}^* \rangle \\
&  = - \langle\eta^{1,*}, u - u_{\max} \rangle + 0 \geq 0\,.
\end{align*}
The same holds for $\eta_{(n)}^{*,2}$, leading to 
\begin{equation*}
L_n(u) - L_n(u_{(n)}^*) \geq  \left \langle u - u_{(n)}^*, \left(F^\top_cF_c + \nu I_\d \right)(u - u_{(n)}^*) \right \rangle\,.
\end{equation*}
\end{proof}

We introduce a variance--like measure $\V^*_{(n)}: \mathcal{P}(\Rd) \to [0,\infty]$ at the solution $u^*_{(n)}$ of the problem at the step $n$ by
\begin{equation}
\label{eq:pseudoV}
\V^*_{(n)}[f]: = \int |u  - u^*_{(n)}|\,f(du)\,.
\end{equation}

Similarly to \cite{fornasier2021consensusbased}, we derive a quantitative version of the Laplace principle, but tailor it to functions satisfying \eqref{eq:grow2}.

\begin{proposition}[Quantitative Laplace principle]
\label{prop:laplace_new}
Let $f \in \mathcal{P}(\Rd)$ with $u_{(n)}^* \in \supp(f)$ and fix $\alpha >0$. If $L_n$ satisfies \eqref{eq:grow2}, then for any $\ve>0$ it holds
\[
|m_{L_n}^\alpha [f] - u_{(n)}^*| \leq \frac {\ve}2  + \frac{\exp(-\alpha \lambda_A^{\min}(\ve/4)^2)}{f(B_{R^n_{\ve}}(u_{(n)}^*))}\, \V_n^*[f]\,
\]
with 
\begin{equation}
R_{\ve}^n \leq \min \left\{\ve^2 \frac{\lambda_A^{\min}}{8(\lambda_A^{\max} + |\eta_{(n)}^{1,*}| + |\eta_{(n)}^{2,*}|)}  \;,\; 1 \right\}\,. 
\label{eq:Rchoice}
\end{equation}
\end{proposition}
\begin{proof}
We note that the Boltzmann-Gibbs distribution is invariant with respect to translation of the objective function. Therefore, it holds
\[ m_{L_n}^\alpha[f] = \frac{\int u e^{-\alpha (L_n(u) - L_n(u_{(n)}^*))} \,f(du)}{\int e^{-\alpha (L_n(u) - L_n(u_{(n)}^*))} \,f(du)} \] 
for any $f \in \mathcal{P}(\Rd)$. Consider two radii $R_1> R_2>0$. By Jensen's inequality, we have
\begin{align*}
|m_{L_n}^\alpha[f] - u_{(n)}^*| &\leq \int |u - u^*_{(n)}| \frac{e^{-\alpha (L_n(u) - L_n(u_{(n)}^*))}}{\int e^{-\alpha (L_n(u) - L_n(u_{(n)}^*))}f(du)}f(du)  \\
& \leq R_1 + \int_{B_{R_1}^c(u_{(n)}^*)} |u - u^*_{(n)}| \frac{e^{-\alpha (L_n(u) - L_n(u_{(n)}^*))}}{\int e^{-\alpha (L_n(u) - L_n(u_{(n)}^*))}f(du)}f(du)\,,
\end{align*}
where we used that 
\[\int_{B_{R_1}(u_{(n)}^*)} |u - u^*_{(n)}| \frac{e^{-\alpha (L_n(u) - L_n(u_{(n)}^*))}}{\int e^{-\alpha (L_n(u) - L_n(u_{(n)}^*))}f(du)}f(du) \leq R_1\,. \] 
To estimate the normalizing constant of Boltzmann-Gibbs distribution, we apply the growth estimates \eqref{eq:grow2}, together with Markov's inequality. Let $\eta_{(n)}:= |\eta_{(n)}^{1,*}| + |\eta_{(n)}^{2,*}| $, this leads to 
\begin{align*}
\int e^{-\alpha (L_n(u) - L_n(u_{(n)}^*))}\,f(du)  &\geq \int e^{-\alpha \lambda^{\max{}}_A | u - u_{(n)}^*|^2 - \alpha\eta_{(n)}|u - u_{(n)}^*|} f(du) \\
& \geq e^{-\alpha \lambda^{\max{}}_A  R_2^2  - \alpha \eta_{(n)} R_2} f\left(B_{R_2}(u_{(n)}^*)\right)\,.
\end{align*}
By exploiting the lower bound in \eqref{eq:grow2}, we obtain 
\begin{align*}
 \int_{B_{R_1}^c(u_{(n)}^*)} & |u - u^*_{(n)}| \frac{e^{-\alpha (L_n(u) - L_n(u_{(n)}^*))}}{\int e^{-\alpha (L_n(u) - L_n(u_{(n)}^*))}\,f(du)}\,f(du) \\
 & \leq
 \frac{e^{\alpha \lambda_A^{\max{}}R_2^2 +  \alpha\eta_{(n)}R_2}}{f\left(B_{R_2}(u_{(n)}^*) \right)}
 \int_{B_{R_1}^c(u_{(n)}^*)} |u - u^*_{(n)}| e^{-\alpha (L_n(u) - L_n(u_{(n)}^*))}\,f(du) \\
 & \leq  \frac{e^{\alpha (\lambda_A^{\max{}}R_2^2 + \alpha\eta_{(n)}R_2 - \lambda_A^{\min{}}R_1})}{f\left(B_{R_2}(u_{(n)}^*) \right)}
 \int_\Rd |u - u^*_{(n)}|\,f(du)\,.
\end{align*} 
Together, we have derived the upper bound 
\[ |m_{L_n}^\alpha[f] - u_{(n)}^*| \leq  R_1 + \frac{e^{\alpha (\lambda_A^{\max{}}R_2^2 + \alpha \eta_{(n)} R_2  - \lambda_A^{\min{}}R^2_1})}{f\left(B_{R_2}(u_{(n)}^*) \right)}
 \V^*_{(n)}[f]\,.\] 
The desired estimate is obtained by the choices $R_1 = \ve/2$ and $R_2 =R^n_{\ve},$ respectively.
\end{proof}

    To apply Proposition~\eqref{prop:laplace_new} for all $f_{(n,k)}$, $n\in \mathbb{N}_0$, $k\in [\overline{k}]:= \{0, \dots, \overline{k} \}$, we  bound from below the quantity $f_{(n,k)}(B_r(u_{(n)}^*))$, for a possibly small radius $r>0$.

\begin{lemma}(Lower bound on mass around $u_{(n)}^*$)
\label{l:mass}
Assume the random noise $\theta$ is Gaussian, $\mu = \mathcal{N}(0, I_\d)$, and consider a fixed radius $r>0$. If $\supp(f_{(0,0)}) = \mathcal{U}$, then there exists a positive constant $\delta_r =\delta_r (\diam(\mathcal{U}),\sigma,\d,f_{(0,0)})$, such that for all $n\in \mathbb{N}_0$, $k\in [\overline{k}]$ it holds
\begin{equation} \notag
f_{(n,k)}\left(B_r(u_{(n)}^*) \right) \geq \delta_r \,. 
\end{equation}
\end{lemma}

\begin{proof} 
Consider $n,k \in \mathbb{N}$ and let $\chi_{D}$ denote the indicator function for the set $D$. By the weak formulation \eqref{eq:MFdiscr} with $\phi(\cdot) = \chi_{B_r(u_{(n)}^*)}(\cdot)$, we get 
\begin{equation}
f_{(n,k)}(B_r(u^*_{(n)})) = \iint \chi_{B_r(u^*_{(n)})} (u')\, f_{(n,k-1)}(du)\, \mu(d\theta)\,. 
\end{equation} 
We introduce the set 
\[ \Theta_{u} := \{\theta \in \Rd\;:\;u' \in B_r(u^*_{(n)})   \}\,.
\] 
Thanks to the convexity of $\mathcal{U}$, we have that $\proj_{\mathcal{U}}$ is non-expansive and $|\proj_{\mathcal{U}}(v) - u^*_{(n)}| \leq |v - u^*_{(n)}|$.
By definition of $u'$ (see \eqref{eq:collCBO}) condition $u' \in B_r(u^*_{(n)})$ is then less restrictive than 
\[  \sigma \sqrt{\tau} \theta \in B_r \left( u^*_{(n)} -   u + \lambda \tau(m_{L_n}^\alpha[f_{(n,k)}] - u) \right) =: B_r( \tilde u )\,. \]
Since $\mu = \mathcal{N}(0, I_\d)$,
\begin{align*}
\mu(\Theta_{u}) &\geq \frac{1}{\sqrt{2\pi \d}}  \int \chi_{B_r(\tilde u)}(\sigma\sqrt{\tau} \theta) e^{- |\theta|^2/2} d\theta \\
& \geq  \frac{e^{- \diam(\mathcal{U})^2/2}}{\sqrt{2\pi \d }}  \int \chi_{B_{r/(\sigma \sqrt{\tau})}(\tilde x/\sigma)}(\theta) d\theta\\
& \geq \frac{e^{- \diam(\mathcal{U})^2/2}}{\sqrt{2\pi \d }} \mathcal{L}^\d\left( B_{r/(\sigma \sqrt{\tau})}(0)\right):=\delta_{r,1}
\end{align*}
where $\mathcal{L}^\d$ is the $\d$-dimensional Lebesgue measure. Next, we apply Fubini's Theorem to conclude that for all $k\in [\overline k]$
\begin{align*}
f_{(n,k)}(B_r(u^*_{(n)})) &= \iint \chi_{B_r(u^*_{(n)})} (u') \,f_{(n,k-1)}(du)\,\mu(d\theta) \\
&= \int \mu(\Theta_u) \,f_{(n,k-1)}(du) \\
&\geq\int \delta_r \,f_{(n,k-1)}(du)  = \delta_{r,1} \,. 
\end{align*} 
For the case $n \in \mathbb{N}, k = 0$, we note that $f_{(n,0)} = f_{(n-1, \overline{k})}$ due to \eqref{eq:CBOinitial}, and so the above lower bound still holds. For the case $(n,k) = (0,0)$, thanks to the compactness of $\mathcal{U}$ and $\supp(f_{(0,0)}) = \mathcal{U}$, we note that $f_{(0,0)}(B_r(\hat u))> \delta_{r,2}$ for all $\hat u \in \mathcal{U}$, for some $\delta_{r,2} >0$. Therefore, the choice $\delta_r :=\min \{\delta_{r,1},\delta_{r,2}\}$ leads to the desired lower bound.
\end{proof}

Next, the evolution of $\V_{(n)}^*[f_{(n,k)}]$ as $k$ increases, is studied. 

\begin{lemma} \label{l:Vdecay}
Assume $|m^\alpha_{L_n}[f_{(n,k)}] - u_{(n)}^*| \leq B$ for all $k \in [\overline{k}]$, for some $B \geq 0$. Then,
\[ V^*_{(n)}[f_{(n,\overline{k})}] \leq e^{- \lambda \overline{k} \tau}\V^*_{(n)}[f_{(n,0)}] + \left(B + \frac{\sigma}{\lambda \sqrt{\tau}}\right) (1 - e^{-\lambda \overline{k} \tau})\,. \] 
\end{lemma}

\begin{proof}

We apply \eqref{eq:MFdiscr} with $\phi(u) = |u|$ and obtain, since the projection is non-expansive,
\begin{align*}
V_{(n)}^*[f_{(n,k+1)}] &= \iint |u'- u_{(n)}^*| \,f_{(n,k)}(du)\,\mu(d\theta) \\
&  \leq \iint |u + \lambda \tau (m^\alpha_{L_n}[f_{(n,k)}] - u) + \sigma \sqrt{\tau}\theta - u_{(n)}^*| \,f_{(n,k)}(du)\,\mu(d\theta) \\
&  \leq (1-\lambda \tau)\V^*_{(n)}[f_{(n,k)}] + \lambda \tau |m^\alpha_{L_n}[f_{(n,k)}] - u_{(n)}^*| + \sigma \sqrt{\tau} \int |\theta|\, \mu(d\theta)\,.
\end{align*}
Using $\mu = \mathcal{N}(0,I_\d)$, and $|m^\alpha_{L_n}[f_{(n,k)}] - u_{(n)}^*|\leq B$, an application of the discrete Grönwall inequality leads to the desired estimate.
\end{proof}

Next, we show that $m_{L_n}^\alpha[f_{(n,k)}]$ is arbitrary close to the minimizer $u^{*}_{(n)}$ provided that the parameter $\alpha$ of the Gibbs-Boltzmann distribution is sufficiently large.

\begin{theorem}
\label{t:malpha}
Consider $f_{(n,k)} = \law(\overline{U}_{(n,k)})$ updated by \eqref{eq:discr_mono}, with $f_{(n,0)} = f_{(n-1,\overline{k})}$ for all $n \in \mathbb{N}$, and $\supp(f_{(0,0)}) = \mathcal{U}$. Let Assumption \ref{asm:linearu} hold, and the MPC step $n\in \mathbb{N}_0$ be fixed. 
For any accuracy $\ve>0$, there exists a value $\alpha^n_0 = \alpha^n_0(\ve,R^n_\ve, \diam(H), \tau, \lambda, \sigma, \d)$ such that for all $\alpha \geq \alpha^n_0$ it holds 
\begin{equation*}
| m^\alpha_{L_n}[f_{(n,k)}] - u_{(n)}^*| \leq \ve \, , 
\end{equation*} 
for all $k \in [\overline{k}]$. Moreover, it holds
\begin{equation}
\lim_{\overline{k} \to \infty} \V_{(n)}^*[f_{(n,\overline{k})}]  \leq  \ve + \frac{\sigma}{\lambda \sqrt{\tau}}\,. 
\label{eq:Vlim}
\end{equation} 
\end{theorem}

\begin{proof} We prove the statement by induction on $k$. Assume that $|m^\alpha_{L_n}[f_{(n,h)}] - u_{(n)}^*| \leq \ve$ for all $0\leq h \leq k-1$. By Lemma \eqref{l:Vdecay} with $B = \ve$, it holds 
\[
\V_{(n)}^*[f_{(n,k)}] \leq \max\left \{\V_{(n)}^*[f_{(n,0)}],\, \ve + \frac{\sigma}{\lambda \sqrt{\tau}}\right\}.
\]

By applying Proposition \ref{prop:laplace_new} (quantitative Laplace principle) and the estimate provided by Lemma \ref{l:mass},  we obtain
\begin{align}
|m^\alpha_{L_n} [f_{(n,k)}] - u_{(n)}^*| &\leq \frac {\ve}2  + \frac{\exp(-\alpha \lambda_A^{\min}(\ve/4)^2)}{\delta_{R^n_\ve}} \V^*_{(n)}[f_{(n,k)}] \notag\\
&\leq \frac {\ve}2  + \frac{\exp(-\alpha \lambda_A^{\min}(\ve/4)^2)}{\delta_{R^n_\ve}}\max\left \{\V_{(n)}^*[f_{(n,0)}],\, \ve + \frac{\sigma}{\lambda \sqrt{\tau}}\right\}  \label{eq:uppbound2}\,.
\end{align}
We note that the second term can be made  arbitrary small by taking $\alpha$ large. In particular, we have
\[|m^\alpha_{L_n} [f_{(n,k)}] - u_{(n)}^*| \leq \ve\,\]
provided $\alpha \geq \alpha^n_0$ with
\begin{equation}
 \alpha^n_0 
 =  \frac{16}{\ve^2 \lambda_A^{\min}} \log \left( \frac2{\ve \delta_{R^n_\ve}} \max\left \{\V_{(n)}^*[f_{(n,0)}],\, \ve + \frac{\sigma}{\lambda \sqrt{\tau}}\right\}\right) \,.
 \label{eq:alpha0n}
\end{equation}
Since \eqref{eq:uppbound2} holds for $k = 0$ too, we have proved that $|m^\alpha_{L_n} [f_{(n,k)}] - u_{(n)}^*| \leq \ve$ for all $k \in [\overline{k}]$. The additional estimate \eqref{eq:Vlim} follows directly  by Lemma \ref{l:Vdecay}, since $\alpha_0^n$ is independent on $\overline{k}$.
\end{proof}

\subsection{Uniform convergence for all MPC sub-problems}

In the estimates derived in Theorem \ref{t:malpha}, we note that the choice of the parameter $\alpha \geq \alpha_0$ depends on the MPC step $n$ through the  two terms: $R_\ve^n$, and $\V_{(n)}^*[f_{(n,0)}]$.  We will show in the following that both quantities are independent of $n.$ In order to obtain 
those estimates, we impose further assumptions on the state space $\mathcal{X}$, and the number $\overline{k}$ of CBO iterations per MPC problem is sufficiently large.

\begin{lemma} \label{l:etas}
Assume the state space $\mathcal{X}$ to be bounded, 
then there exists  a radius $R_\ve>0$ independent on $n$, such that  \eqref{eq:Rchoice} holds for all $n \in \mathbb{N}_0$.
\end{lemma}
\begin{proof}
We note that the choice \eqref{eq:Rchoice} of $R_\ve^n$ depend on $n$ through the Lagrange multipliers $\eta_{(n)}^{1,*}, \eta_{(n)}^{2,*}$, and, in particular, by $|\eta_{(n)}^{1,*}| + |\eta_{(n)}^{2,*}|$. First of all, we note that, since $u_{\min}< u_{\max}$, only one constraint can be active in each component, leading to $\langle\eta_{(n)}^{1,*}, \eta_{(n)}^{2,*}\rangle =0.$ As a consequence, $|\eta_{(n)}^{1,*}|^2 + |\eta_{(n)}^{1,*}|^2 = |\eta_{(n)}^{1,*} - \eta_{(n)}^{2,*}|^2$. Next, we  bound their difference. This  is obtained by the  KKT condition \eqref{eq:kkt}
\begin{align*}
 | \eta_{(n)}^{*,1} - \eta_{(n)}^{*,2}| &= \left|2\left( F_c^\top F_c + \nu I_\d \right) u_{(n)}^* - 2F_c^\top \left( \Phi_s (x_{(n)})  -  x^{\refup}(t_{n+1})  \right) \right|\\
 & \leq 2 \lambda_A^{\max} |u_{(n)}^*| + 2\lambda_{F_c^\top}^{\max}|\Phi_s (x_{(n)})  -  x^{\refup}(t_{n+1})| \\
 & \leq 2 \lambda_A^{\max} \diam(\mathcal{U}) + 2\lambda_{F_c^\top}^{\max} \diam (\mathcal X)\,.
\end{align*}
We can therefore obtain an estimate for $R_\ve^n$  independent on $n$.
\end{proof}

\begin{theorem} Under the assumptions to Theorem \ref{t:malpha}, let additionally assume $\mathcal{X}$ to be bounded. Fix an arbitrary accuracy $\ve$ and let $R_\ve$ be given by Lemma \ref{l:etas}. Then, if $\alpha> \alpha_0$ 
with 
\begin{equation}
\label{eq:alpha0}
\alpha_0  =  \frac{16}{\ve^2 \lambda_A^{\min}} \log \left( \frac2{\ve \delta_{R_\ve}} \left (\diam(\mathcal{U})+ \ve + \frac{\sigma}{\lambda \sqrt{\tau}}\right)\right) 
\end{equation}
and
\begin{equation} \label{eq:kchoice}
\overline{k} \geq \frac1{\lambda \tau} \log \left(\frac{\diam(\mathcal{U})}\ve \right)\,,
\end{equation}
it holds for all $n \in \mathbb{N}_0$
\begin{align*}
|u_{(n)} - u_{(n)}^*|& \leq \ve \\
\V^*_{(n)}[f_{(n, \overline{k})}] &\leq 2\ve + \frac{\sigma}{\lambda \sqrt{\tau}}\,.
\end{align*}
\end{theorem}

\begin{proof}
We note that $\V^*_{(n)}[f] \leq \diam(\mathcal{U})$ for any $f\in \mathcal{U}$ and $n\in\mathbb{N}_0$. Together with $R_\ve \leq R_\ve^n$, this leads to $\alpha_0 \geq \alpha_0^n$ for all $n$.

Therefore, we apply Theorem \ref{t:malpha} (recall $u_{(n)}  = m^\alpha_{L_n}[f_{(n,\overline{k})}]$ due to Assumption \ref{asm:linearu}) and, in turn, Lemma \ref{l:Vdecay} with $B = \ve$. This leads to the bound
\[ \V^*_{(n)}[f_{(n,\overline{k})}] \leq e^{- \lambda \overline{k} \tau}\V^*_{(n)}[f_{(n,0)}] + \left(\ve + \frac{\sigma}{\lambda \sqrt{\tau}}\right) (1 - e^{-\lambda \overline{k} \tau})
\,.
\]
The desired upper bound is obtained by using again $\V^*_{(n)}[f] \leq \diam(\mathcal{U})$, and the choice of $\overline{k}$ given by \eqref{eq:kchoice} which shows that the first term is  smaller than $\ve$.
\end{proof}

\section{Numerical experiments}
\label{sec:num}

We validate the CBO solver for MPC, by applying Algorithm \ref{alg:mpc-cbo} to the control of a simulated Continuous Stirred-Tank Reactor (CSTR) nonlinear system \cite{biagiola2005use,stahl2011}.

\subsection{Description of CSTR}
The system description is given in \cite{stahl2011}. The CSTR is a chemical reactor where a chemical $A$ is supplied at a concentration $C_f$ and temperature $T_0$. An irreversible, first-order exothermic reaction $A \to B$, for some chemical $B$ takes place. Being the reactor stirred, the concentration $C$, as well as the temperature $T$ is assumed to be uniform. The system is controlled in order to keep the concentration $C$ at a given value $C^\refup$. This is formulated  as controlling a coolant stream $q_c$ flowing around the CSTR.

The system evolution is given by 
\begin{equation}\label{eq:cstr}
\begin{split}
\frac{dC}{dt} & = \frac{q}{V} \left( C_f - C \right) - k_0 C e^{-E/(RT)}
\\
\frac{dT}{dt} & = \frac{q}{V}\left( T_0 - T\right) - \frac{k_0 \Delta H}{ \rho c_p} C e^{-E/(RT)} + \frac{\rho_c c_{pc} q_c}{\rho c_p V} \left(1 - \exp\left(-\frac{hA}{q_c \rho_c c_{pc}} \right)     \right) (T_{C_0} - T)\,.
\end{split}
\end{equation}
We refer to Table \ref{table:cstr} for the complete list of CSTR parameters. The reference state is given by $C^\refup(t) = 0.1 $ for $t \in [0,3]$ and $C^\refup(t) =  0.12, $ for $t \in [3,6.5]$ with corresponding reference coolant steam $q_c^\refup = 103.411, 108.1.$

\begin{table}
\centering
\begin{tabular}{|l l l l |} 
 \hline
 Parameter & Description & Unit & Value \\ [0.5ex] 
 \hline
 $C$              & Concentration            & mol/l & -- \\ 
 $C_f$           & Feed concentration    & mol/l & 1 \\
 $T$         	     & Temperature               & K     & -- \\
 $T_0$    	     & Feed temperature      & K       & 350 \\
 $q_{c}$         & Coolant stream temp.  & K     & --    \\
 $V$              & Reactor volume          & l       & 100\\ 
 $q$               & Process flow rate        & 1/min    &100    \\
 $T_{C_0}$    & Inlet coolant temp. & K    & 350   \\
 $hA$            &  Heat transfer term         &Kcal/min   & $7 \times 10^5$   \\
 $k_0$           & Reaction rate         & 1/ min     & $7.2\times10^{10}$   \\
 $E/R$           & Activation energy  & K     & $10^4$   \\
 $\Delta H$     & Heat of reaction          &  cal/mol   &  $-2\times 10^5$   \\
 $\rho, \rho_c$& Liquid densities          &   g/l  &  $10^3$  \\
 $c_p, c_{pc}$&   Specific heats        & k cal/g    & 1   \\
 [1ex] 
 \hline
\end{tabular}
\caption{Parameters used for simulation of CSTR.}
\label{table:cstr}
\end{table}

\subsection{Results}

To apply the control strategy, we assume the state variables $C_{(m)} = C(t_m),T_{(m)}  = T(t_m)$ to be measured
at a sampling time $\Delta t= 0.05$ minutes. We consider a time window of $\overline n \Delta t = 6.5$ minutes. 
The nonlinear ODE system is solved via explicit Euler method with a times step of $10^{-3}$  minutes.
In the experiments, we use a prediction horizon $n_p = 10$ and a loss function at step $n$ given by 
\[ L_n((q_c)_{n:n+n_p}) =
\sum_{m = n+1}^{n + n_p +1 } |C_{(m)} -  C^{\refup}(t_{m})|^2 + \nu \sum_{m = n}^{n+n_p} |(q_c)_{(m)} - q_c^\refup(t_m)|^2 
\]
with control variable $u = q_c \in [20,200]$ and $\nu = 1$. 

The CBO  system is initialized according to  $U^i_{(0,0)}\sim q_c^\refup(0) + \textup{Unif}([-0.5, 0.5]^{n_p})$ for $i = 1, \dots, N$, and the initial state variables are given by $x_{(0)} = (C_{(0)}, T_{(0)}) = (0.1,438.54)$. We use non-degenerate diffusion for the CBO agents given by
\[D  = D^i_{(n,k)}  = (m^\alpha_{L_n}[f_{(n,k)}] - U_{(n,k)}^i) + \tilde\sigma \mathbf{1}_\d\]
for $\hat \sigma = 10^{-3}.$ 
The above choices correspond to a strategy between a pure white noise, as considered in analyzed in Section \ref{sec:analysis}, and the anisotropic noise proposed in \cite{carrillo2019consensus}.
The remaining CBO parameters are set to $\alpha = 10^5, \lambda = 1, \sigma = 3, \tau = 0.1.$

\begin{figure}
\centering
\includegraphics[width = 9 cm]{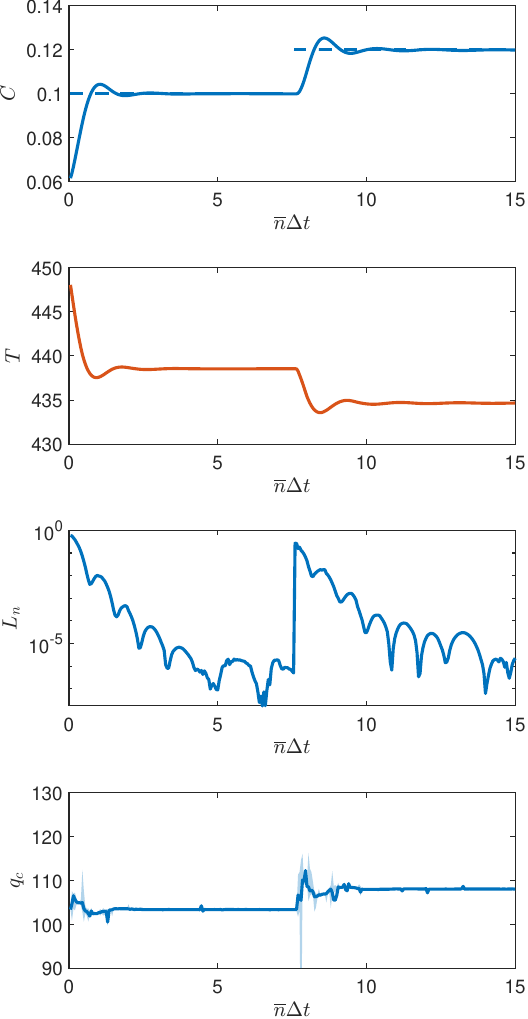}
\caption{Results for a single run of Algorithm \ref{alg:mpc-cbo}. $N = 32$ agents are used, with $\overline{k} = 10$ and MPC horizon $n_p = 10$.}
\label{fig:singlerun}
\end{figure}

For all  experiments we set a MPC control horizon to $n_p = 10$.
Figure \ref{fig:singlerun} shows the results of a single experiment with $N = 32$ CBO agents, and $\overline k = 10$ CBO iterations per sub-problem. Considering the  evolution of $L_n$, we note that Algorithm \ref{alg:mpc-cbo} is able to control the system and to reach a MPC loss accuracy of order $10^{-7}$. 

\begin{figure}
\centering
\includegraphics[width = 9cm]{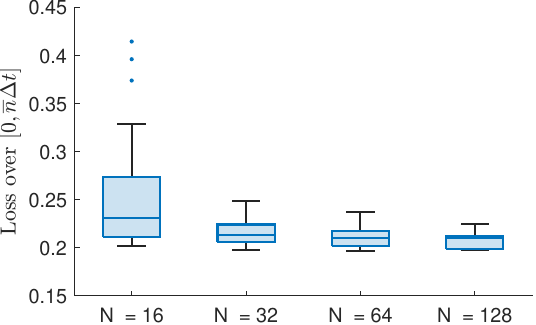}
\caption{Statistics of the of the algorithm's perfomance for different numbers $N$ of agents employed, computed over 20 experiments. $\overline{k} = 15$ CBO iterations are performed per MPC sub-problem.}
\label{fig:N}
\end{figure}

As the proposed theoretical analysis relies on a mean-field approximation of the CBO agent system, we investigate how the number of agents $N$ affects the performance of the algorithm. As shown in Figure \ref{fig:N}, employing more than $32$ agents leads, on average, to a marginal improvement of the total loss $L_n$. On the other hand, the inter-quantile ranges diminish for large $N$, showing that the algorithm performance is  consistent provided  more agents are used. This is coherent with the propagation of chaos assumption, which implies  that the CBO  distribution $f_{(n,k)}$ is deterministic for $N \to \infty$.

\begin{figure}
\centering
\includegraphics[width = 9cm]{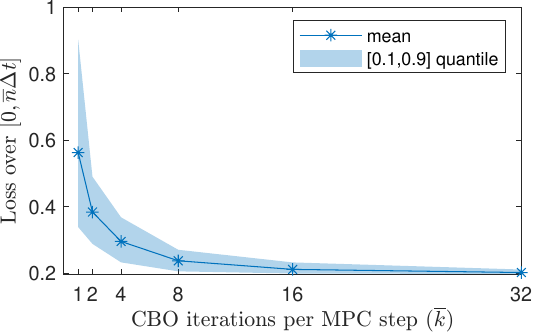}
\caption{Statistics of the of the algorithm's performance for different numbers $\overline{k}$ of CBO iterations per MPS problems, computed over 30 experiments. $N = 32$ agents are used.}
\label{fig:maxit}
\end{figure}

Figure \ref{fig:maxit} shows that  the number CBO iterations $\overline{k}$ per MPC sub-problem impacts the performance of the proposed strategy. Using more than $16$ iterations leads, on average, to only a marginal improvement in total loss $L_n$ achieved. The initialization  \eqref{eq:CBOinitial}, that is, $f_{(n+1,0)} = f_{(n, \overline{k})}$, allows the CBO agent system to sequentially solve  many MPC subproblem with only a few iterations. For a general  reference to the  CBO literature, we refer to \cite[Section 4.2]{carrillo2019consensus}, where $10^4$ steps are used for $20$-dimensional problems. We also note that the inter-quantile range decreases as $\overline{k}$ increases, showing, again, a trade-off  between computational cost and fluctuations of the algorithm performance.

\section{Conclusion}

In this work, we  coupled a CBO agent system with the MPC strategy to derive a novel method for the online control of time-discrete nonlinear systems. The final agents' locations are used as  initialization of the agent system for the subsequent MPC sub-problem. Intuitively, this means that the CBO system runs parallel to the state dynamics, but using more frequent updates.

A theoretical analysis of systems with linear additive control, shows that the convergence properties of the CBO interaction do not deteriorate as the state system evolves. The analysis is performed by relying on a mean-field approximation of the multi-agent system. We are  able to bound the difference between the applied and optimal control.

Numerical test for the CSTR shows the validity of the proposed algorithm, that is able to control the system with only $N = 32$ agents and a model predictive horizon of $n_p = 10$ time steps. The initialization of the agent system after every sub-problem, shows that few CBO iterations are sufficient. 

The proposed coupled state-control system could also be suitable for  time-continuous systems. In future work, we aim to derive a time-continuous MPC strategy independent of the sampling period.


\bibliographystyle{elsarticle-num} 
\bibliography{bibfile}


%
%
%
%


\end{document}